\documentclass[conference]{IEEEtran}
\ifCLASSINFOpdf
   \usepackage[pdftex]{graphicx}
\else
\fi
\usepackage[cmex10]{amsmath}
\usepackage{amsfonts, amssymb}
\usepackage{color}
\usepackage{bbm}
\usepackage{amsthm}

\graphicspath{{fig/}}

\makeatletter
\newcommand*{\ov}[1]{\m@th\overline{\mbox{#1}\raisebox{.8em}{}}}
\begin{document}

\newtheorem{lem}{Lemma}
\newtheorem*{mydef}{Definition}
\newtheorem{thm}{Theorem}
\newtheorem{prop}{Proposition}
\newtheorem*{notat}{Notation}

\abovedisplayskip=2pt
\belowdisplayskip=2pt
\abovedisplayshortskip=2pt
\belowdisplayshortskip=2pt

%
\title{Construction of Orthonormal Quasi-Shearlets based on quincunx dilation subsampling}

\author{\IEEEauthorblockN{Rujie Yin}
\IEEEauthorblockA{Department of Mathematics, Duke University, USA\\
Email: rujie.yin@duke.edu,
}
}


%


\IEEEoverridecommandlockouts
\IEEEpubid{\makebox[\columnwidth]{978-1-4799-7492-4/15/\$31.00~
\copyright2015
IEEE \hfill  accepted version of \cite{sampta15-shearlet}} \hspace{\columnsep}\makebox[\columnwidth]{ } }

\maketitle

\begin{abstract}
We consider the construction of orthonormal directional wavelet bases in the multi-resolution analysis (MRA) framework with quincunx dilation downsampling. We show that the Parseval frame property in MRA is equivalent to the identity summation and shift cancellation conditions on $M$ functions, which essentially characterize the scaling (father) function and all directional (mother) wavelets. 

Based on these two conditions, we further derive sufficient conditions for orthonormal bases and build a family of quasi-shearlet orthonormal bases, that has the same frequency support as that of the least redundant shearlet system.  In addition, we study the limitation of our proposed bases design due to the shift cancellation conditions.
 \end{abstract}


%
\IEEEpeerreviewmaketitle

\section{Introduction}
In image compression and analysis, 2D tensor wavelet schemes are widely used, despite having the severe drawback of poor orientation selectivity: even for ``cartoon images'' (in which the image is piecewise smooth, with jumps occurring possibly along piecewise $C^2$-curves), only horizontal or vertical edges can be well represented by tensor wavelets. To resolve this problem, several different constructions of directional wavelet systems have been proposed in recent years, with promising performance. 
Shearlet \cite{shearlet12book,easley2008sparse} and curvelet \cite{candes2006fast} systems construct a multi-resolution partition of the frequency domain by applying shear or ratation operators to a generator function in each level; contourlets \cite{do2005contourlet} combine the Laplacian pyramid scheme with the directional filter approach; dual-tree wavelets \cite{selesnick2005dual} are linear combination of 2D tensor wavelets (corresponding to multi-resolution systems) that constitute an approximate Hilbert transform pair. 

Among these,  shearlet and curvelet systems have optimal asymptotic rate of approximation for ``cartoon images'' \cite{guo2007optimally,candes2005curvelet} and have been successfully applied to image denoising \cite{easley2009shearlet}, image restoration \cite{candes2002new} and image separation \cite{kutyniok2012image}. Despite these striking theoretical potential for curvelets and shearlets, their implicit (often high) redundancy impedes their practical usage. 
Depending on the property of the generator function and the number of directions, available shearlet and curvelet implementations have redundancy of at least 4; in addition, the factor may grow exponentially in the number of decomposition levels.

Although redundancy is exploited in image processing tasks such as denoising, restoration and reconstruction, a basis decomposition like the 2D tensor wavelets is preferred in image compression and other tasks where computation cost is of concern.

Our goal is to construct a directional wavelet system that has similar orientation selectitvity to the shearlet system in the first ``dyadic ring'' in the frequency domain, but with much less and ideally no redundancy. Thus, it is reasonable to call such wavelet schemes ``quasi-shearlets''. The paper is organized as follows: we first set up the framework of a dyadic MRA with quincunx dilation downsampling and show that such non-redundant downsampling is critical for the partition of frequency domain to be similar to a shearlet system. Then, we derive two essential conditions, {\it identity summation} and {\it shift cancellation}, for perfect reconstruction in this MRA with critical downsampling. These conditions lead to the classification of {\it regular/singular} boundaries of the frequency partition and the corresponding smoothing techniques to improve spatial localization of the quasi-shearlets. Finally, we present a family of quasi-shearlet orthonormal bases constructed in the introduced framework.

\section{Framework setup}\label{sec: framework}
In this section, we summarize the notions of MRA associated to a 2D wavelet system, the matrix representation of a sublattice of the integer grid $\mathbb{Z}^2$ and the relation between frequency domain partition and sublattice with critical downsampling.

\subsection{Multi-resolution analysis and critical downsampling}
%

We start by assuming that $\phi$ is a scaling function for a MRA in $L^2(\mathbb{R}^2)$; more precisely, $\Vert\phi\Vert_2=1$, and there exists a $2\times 2$ matrix $D$ with integer entries s.t. the rescaled $\phi_1(\boldsymbol{x}) = |D|^{-1/2}\phi(D^{-1}\boldsymbol{x})$ can be written as a linear combination of shifted $\phi_{0,\boldsymbol{k}} = \phi(\boldsymbol{x}-\boldsymbol{k}), \, \boldsymbol{k}\in \mathbb{Z}^2$, where $|D|$ is the determinant of $D$. It then follows that
\begin{align}\label{eq: m0}
\hat{\phi}(D^T\boldsymbol{\xi}) = M_0(\boldsymbol{\xi})\hat{\phi}(\boldsymbol{\xi}),
\end{align}
where $M_0(\boldsymbol{\xi}) = M_0(\xi_1,\xi_2)$ is $2\pi-$periodic in $\xi_1,\xi_2$; its Fourier coefficients are given (up to normalization) by those in the linear combination of $\phi(\boldsymbol{x}-\boldsymbol{k})$ (see above). Assuming that $\phi\in L^1(\mathbb{R})$, so that $\hat{\phi}$ is continuous, and $\int \phi(\boldsymbol{x})d\boldsymbol{x} = 1$, iterating \eqref{eq: m0} shows that $\hat{\phi}$ is uniquely determined as an infinite product of rescaled $M_0$,
$\hat{\phi}(\boldsymbol{\xi}) = (2\pi)^{-1}\prod_{k=1}^{\infty}M_0(D^{-k} \boldsymbol{\xi}).$


The MRA then uses the nested approximation spaces $V_l$ defined by $V_l = \overline{span\{\phi(D^l\boldsymbol{x}-\boldsymbol{k});\boldsymbol{k}\in\mathbb{Z}^2\}}$. Next, we assume that we have wavelet functions $\psi^j$, {\small $1 \leq j \leq J$} together with an integer scaling matrix $\tilde{D}$ and {\small $(2\pi,2\pi)$}-periodic functions $M_j$, satisfying 
\begin{align}\label{eq: mj}
\hspace*{0em}\hat{\psi}^j_1(\boldsymbol{\xi})\propto\,\hat{\psi}^j(\tilde{D}^T\boldsymbol{\xi}) = M_j(\boldsymbol{\xi})\hat{\phi}(\boldsymbol{\xi}),\hspace{1cm} 1\leq j \leq J.
\end{align}
We shall then consider the $\psi_{1,k}^j(\boldsymbol{x}) = \psi_1^j(\boldsymbol{x}-\tilde{D}\boldsymbol{k}),\boldsymbol{k}\in\mathbb{Z}^2,\,1\leq j\leq J,$ and the corresponding spaces $W_1^j = \overline{span\{\psi_{1,k}^j,k\in\mathbb{Z}^2\}}$. Together, the $W_1^j,1\leq j\leq J,$ and $V_1$ span $V_0$. In order to avoid redundancy (or to have {\it critical} downsamping), it is necessary that
$$|D|^{-1} + J|\tilde{D}|^{-1} = 1.$$
Finally, an $L$-level multi-resolution system with base space $V_0$ consists of $\{\phi_{L,k}\,,\psi^j_{l,k}$\,,{\small $\, 1\leq l \leq L,\, k\in \mathbb{Z}^2,\,1\leq j \leq J\}$}, where $span(\phi_{l-1,k}) = span(\phi_{l,k},\psi_{l,k})$. Next, we show a specific choice of the support of $M_j$'s and scaling matrices $D$, $\tilde{D}$ that achieves critical downsampling.

\subsection{Frequency domain partition and sublattice sampling}
Let $\mathbf{S}_0 = [-\pi,\pi)\times[-\pi,\pi)$ be the initial frequency square with respect to the sampling lattice $\mathbb{Z}^2$, which is also where $\hat{\phi}$ is mostly concentrated. Since $\phi_1,\psi^j_1$ and their shifts span the space $V_0$, $supp(\hat{\phi})$ and $supp(\widehat{\psi^j})$, together, should thus cover $\mathbf{S}_0$. Particularly, $\mathbf{S}_0=\bigcup_{0\leq j\leq J} supp(M_j\vert_{\mathbf{S}_0})$, meaning that a partition of $\mathbf{S}_0$ can be generated corresponding to the intersection of $\mathbf{S}_0$ and the concentration area of the $M_j$'s. For simplicity, we abuse the notation $supp(M_j)$ for the concentration area of $M_j\vert_{\mathbf{S}_0}$ hereafter.

\begin{figure}[!t]
\centering
\includegraphics[width=.2\textwidth]{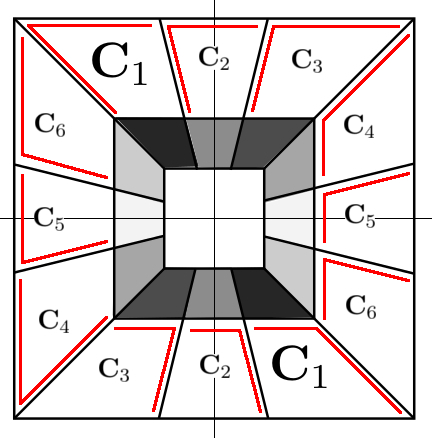}\hspace*{2mm}
\includegraphics[width=.2\textwidth]{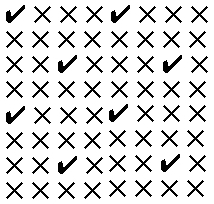}
\caption{Left: partition of $\mathbf{S}_0$ and boundary assignment of $\mathbf{C}_j$'s ( each $\mathbf{C}_j$ has boundaries indicated by red line segments), Right: dilated quincunx sublattice. }
\label{fig: partition}
\end{figure}

To build a quasi-shearlet bases, we choose the partition of $\mathbf{S}_0$ to be that of the least redundant shearlet system, see the left of Fig.\ref{fig: partition}. In this partition, $\mathbf{S}_0$ is divided into a central square $\mathbf{S}_1 = \bigl(\begin{smallmatrix} 2&0\\0&2\end{smallmatrix}\bigr)^{-1}\mathbf{S}_0
$ and a ring, and the ring is further cut into six pairs of directional trapezoids $\mathbf{C}_j$'s by lines passing through the origin with slopes $\pm 1, \pm 3$ and $\pm \frac{1}{3}$. The central square $\mathbf{S}_1$ can be further partitioned in the same way to obtain a two-level multi-resolution system, as shown in the left of Fig.\ref{fig: partition}.

Given this partition, it is straightforward that $supp(M_0) = \mathbf{S}_1$ indicates $D= D_2\doteq\bigl(\begin{smallmatrix} 2&0\\0&2\end{smallmatrix}\bigr)$, i.e. the scaling coefficients are taken on the classical dyadic sublattice. We now claim that each pair of trapezoids $\mathbf{C}_j$ corresponds to $\tilde{D} =Q\doteq \bigl(\begin{smallmatrix} 2&2\\-2&2\end{smallmatrix}\bigr)$, i.e. the shearlet coefficients are taken on the dyadic quincunx sublattice (the right of Fig.\ref{fig: partition}), based on the following lemma.
\begin{mydef}
A compact set $\mathbf{S}$ is a {\it frequency support} of the lattice $\mathcal{L}$ generated by $\boldsymbol{a}_1,\boldsymbol{a}_2$ if $\int_{\mathbf{S}}e^{i\boldsymbol{n}^T\boldsymbol{\xi}}d\boldsymbol{\xi} = \frac{1}{|\boldsymbol{a}_1\times \boldsymbol{a}_2|}(2\pi)^2\delta_{\boldsymbol{n,0}}$.\\
The reciprocal lattice $\tilde{\mathcal{L}}$ of $\mathcal{L}$ is a lattice with basis vectors $b_1,b_2$, s.t. $b_i^Ta_j = 2\pi\delta_{ij}$.
\end{mydef}
\begin{lem}
\label{lem: shifts}
If $A$ is a 2-by-2 matrix with integer entries, then $\mathcal{L}' = A\mathcal{L}$ is a sublattice of $\mathcal{L}$, i.e. $\mathcal{L}'\subset\mathcal{L}$, and there exist shifts $\tilde{\boldsymbol{\lambda}}_m,\,1\leq m\leq |A|-1$, $\tilde{\boldsymbol{\lambda}}_0=0$  s.t.  $\bigcup_{m=0}^{|A|-1} (\mathcal{L}'+\tilde{\boldsymbol{\lambda}}_m) = \mathcal{L}$.\\[.1em]
\end{lem}
\vspace{-2em}
\begin{lem}
\label{lem: sub-shifts}
Given a lattice $\mathcal{L}$, and $\tilde{\boldsymbol{\lambda}}_m,\boldsymbol{\lambda}_m,\,0\leq m\leq |A|-1$, s.t. $\{\tilde{\boldsymbol{\lambda}}_m\}$ and $\vspace{.2em}\{\frac{1}{2\pi}A^T\boldsymbol{\lambda}_m\}$ are sets of shifts associated to $A\mathcal{L}\subset\mathcal{L}$ and $\tilde{\mathcal{L}}\subset(A\mathcal{L})\,\tilde{} = (A^T)^{-1}\tilde{\mathcal{L}}$ respectively, as in Lemma \ref{lem: shifts},
then $\sum_{m}e^{i\boldsymbol{\lambda}_m^T\tilde{\boldsymbol{\lambda}}_n} = 0, \forall 1\leq n\leq |A| -1$.
Furthermore $\mathbf{S}'$ is a frequency support of $\mathcal{L}'$ iff  $\,\bigcup_{m}(\mathbf{S}'+\boldsymbol{\lambda}_m)$ is a frequency support of $\mathcal{L}$.
\end{lem}

According to Lemma \ref{lem: sub-shifts}, $\mathbb{Z}^2$ can be decomposed into a union of shifted sublattices $Q\mathbb{Z}^2+\tilde{\boldsymbol{\lambda}}$; $\mathbf{C}_j$ is the frequency support of $Q\mathbb{Z}^2$ if and only if $\mathbf{S}_0$ is also a union of $\mathbf{C}_j+\boldsymbol{\lambda}$, and the two sets of shifts $\boldsymbol{\lambda}$ and $\tilde{\boldsymbol{\lambda}}$ are linked by $Q$. We check that both decompositions of $\mathbb{Z}^2$ and $\mathbf{S}_0$ hold under matrix $Q$, shifts $\Lambda = \{(\frac{\pi}{2},\frac{\pi}{2}),(\frac{3\pi}{2},\frac{\pi}{2}),(\frac{\pi}{2},\frac{3\pi}{2}),(\frac{3\pi}{2},\frac{3\pi}{2}),$ $(0,0),(0,\pi),(\pi,0),(\pi,\pi)\}$ and each trapezoid pair $\mathbf{C}_j$ for the $j$-th shearlet. Likewise, they hold under matrix $D_2$, shifts $\Gamma = \{(0,0),(0,\pi),(\pi,0),(\pi,\pi)\}$ and $\mathbf{S}_1$ for scaling function. In addition, $|D_2| = 4,\,|Q| = 8$ so that $1/4 + 6\cdot 1/8 = 1$, and the system has critical downsampling.

\section{Parseval frame and bases conditions}
In the rest of the paper, we always consider a multi-resolution system with scaling function $\phi$ and quasi-shearlets $\psi^j$, {\small$ j = 1,\dots,6$} defined by $(M_0, D_2)$ and $(M_j, Q)$, {\small$j = 1,\dots,6$} respectively. Furthermore, the essential support of $M_j$'s corresponds to the partition of $\mathbf{S}_0$ in a shearlet system.

Only the $M$ functions need to be designed to construct this 2D system, which is a close analogue of the classical MRA of 1D wavelets. 
We begin with examining the conditions on $M$ functions such that the system is a Parseval frame.

\subsection{Identity summation and shift cancellation}

In MRA, the Parseval frame property is equivalent to the one-layer perfect reconstruction condition: $\forall f\in L_2(\mathbb{R}^2)$,
$$\textstyle \sum_k\langle f,\phi_{0,k}\rangle\phi_{0,k} = \sum_k\langle f,\phi_{1,k}\rangle\phi_{1,k} + \sum_j\sum_k\langle f,\psi^{j}_{1,k}\rangle\psi^j_{1,k}
.$$
Using \eqref{eq: m0} and \eqref{eq: mj}, this condition on $\phi$ and $\psi^j$'s may be transformed into the condition on $M_j$'s.
\begin{thm}\label{thm: conds}
The perfect reconstruction condition holds iff the following two conditions hold
\begin{align}\label{eq: id-sum}
|M_0(\boldsymbol{\xi})|^2 + \sum_{j = 1}^6|M_j(\boldsymbol{\xi})|^2 = 1
\end{align}
\begin{equation}\label{eq: shift-cancel}
 \begin{cases}
\sum_{j = 0}^6M_j(\boldsymbol{\xi})\overline{M_j(\boldsymbol{\xi} + \boldsymbol{\gamma})} = 0, & \boldsymbol{\gamma}\in \Gamma\setminus\{\boldsymbol{0}\}\\[.5em]
\sum_{j=1}^6M_j(\boldsymbol{\xi})\overline{M_j(\boldsymbol{\xi}+\boldsymbol{\nu})} = 0, & \boldsymbol{\nu}\in\Lambda\setminus\Gamma
\end{cases}
\end{equation}
\end{thm}
In Theorem \ref{thm: conds}, Eq. \eqref{eq: id-sum} is the {\it identity summation} condition, so then the $l_2$ energy is conserved after one-level decomposition; Eq. \eqref{eq: shift-cancel} is the {\it shift cancellation} condition, such that downsampling of scaling and shearlet coefficients is valid. Each $M_j$ contributes a term $M_j(\boldsymbol{\xi})\overline{M_j(\boldsymbol{\xi} + \boldsymbol{\nu})}$ in the cancellation condition of all the shifts corresponding to the downsampling scheme of $M_j$ discussed in section \ref{sec: framework}.

\subsection{Extra condition for basis}
In the 1D wavelet MRA, Cohen \cite{cohen1992biorthogonal} provides a necessary and sufficient condition for a Parseval frame to be an orthonormal basis. The following theorem is a 2D generalization of Cohen's theorem.
\begin{thm}\label{thm: basis cond}
Assume that $M_0$ is a trigonometric polynomial with $M_0(0)=1$, and define $\hat{\phi}(\boldsymbol{\xi})$ recursively by $M_0$.\\
If $\phi(\cdot - n)$ are orthonormal, then $\exists K$ containing a neighborhood of 0, s.t. $\forall\boldsymbol{\xi}\in\mathbf{S}_0,\,\exists \mathbf{n}\in\mathbb{Z}^2, \boldsymbol{\xi}+2\pi\mathbf{n}\in K$ and $\inf_{k>0,\,\boldsymbol{\xi}\in K}|M_0(D_2^{-k}\boldsymbol{\xi})| >0$.  \\
 Furthermore, if $\sum_{\boldsymbol{\gamma}\in \Gamma} |M_0(\boldsymbol{\xi}+\boldsymbol{\gamma})|^2 = 1$, then the inverse is true.
\end{thm}
Because it is difficult to directly design $M_0$ that satisfies the conditions in Theorem \ref{thm: basis cond}, we construct $M$ functions with merely identity summation and shift cancellation constraints and check if the resulting Parseval frame is indeed an orthonormal basis by applying Theorem \ref{thm: basis cond} to $M_0$ afterwards.

\section{M-function design and Boundary regularity}\label{sec: design}
In this section, we first introduce a naive Shannon-type shearlet orthonormal basis, which leads to the critical analysis of the boundary regularity of $\mathbf{S}_1$ and the $\mathbf{C}_j$'s.

\subsection{Shannon-type shearlets}
If each $M_j$ is an indicator function of a piece in the partition of $\mathbf{S}_0$, i.e. $M_0 = \mathbbm{1}_{\mathbf{S}_1},\, M_j = \mathbbm{1}_{\mathbf{C}_j},\,1\leq j \leq 6,\,$ and there is no overlap of two adjacent pieces in the partition on their common boundary,
then the identity summation follows from the nature of partition. Moreover, using the boundary assignment in the left of Fig.\ref{fig: partition}, the shift cancellation holds automatically due to $M_j(\boldsymbol{\xi})\overline{M_j(\boldsymbol{\xi} + \boldsymbol{\nu})}\equiv 0,\,\forall j,\boldsymbol{\nu}\neq 0.$ 
It is also easy to check that the Shannon-type shearlets generated from these $M$ functions form an orthonormal basis.

In addition, one can show, by applying suitable generalizations of Lemmas \ref{lem: shifts}, \ref{lem: sub-shifts} and Theorem \ref{thm: conds} to different dyadic scales, that this Shannon-type construction can be extended to trapezoid-pair supports in the frequency domain by scaling the $\mathbf{C}_1,\dots,\mathbf{C}_6$ dyadically. In particular, aftering scaling by 4, $4\mathbf{C}_j$'s can be further divided in half by a radial line going through the origin; the corresponding subsampling matrix then becomes $QP$, where $P = \bigl(\begin{smallmatrix} 2&0\\0&1\end{smallmatrix}\bigr)$ for horizontal trapezoid pairs, and $P = \bigl(\begin{smallmatrix} 1&0\\0&2\end{smallmatrix}\bigr)$ for vertical trapezoid pairs. The same trick can be repeated every two dyadic scales, leading to an orthonormal shearlet basis in which all shearlets corresponding to sufficient localization in the frequency domain that have disjoint, compact supports.

\begin{figure}[!t]
\centering
\hspace*{-5mm}
\includegraphics[width=.4\textwidth]{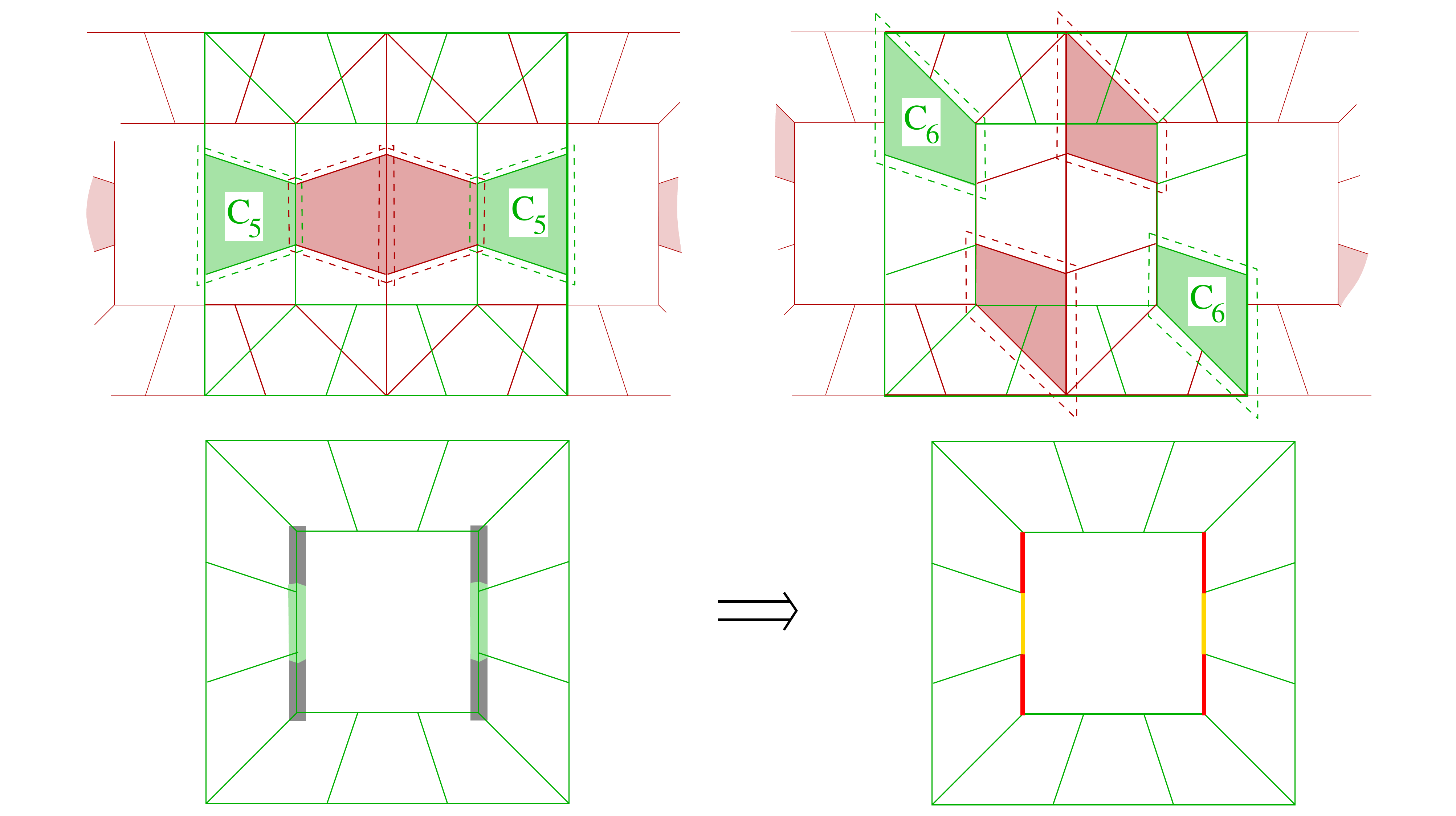}\hspace*{-2mm}
\includegraphics[width=.13\textwidth]{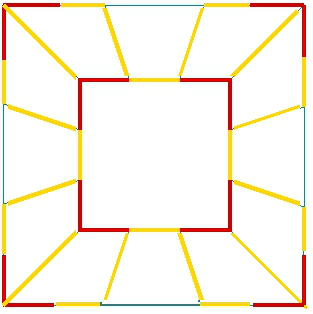}
\caption[caption]{
{\it Left top}: the support of $M_j$ (green) and $M_j(\cdot+(\pi,0))$ (maroon) for $j = 5,6$, after smoothing, overlapping on the vertical boundary at $\xi_1 =  \pm \pi/2$ of $\mathbf{C}_5$ (green) and its shift (red) by $\boldsymbol{\nu} = (\pi,0)$. Note that we have two copies of shifted $\mathbf{S}_0$ (red) overlapping the unshifted $\mathbf{S}_0$(green) due to the $(2\pi,2\pi)$ periodicity of $M_j$. Only $M_0$ and $M_5$ have overlapping smoothed boundaries by this $\boldsymbol{\nu}$. 
\\\hspace{\textwidth}
{\it Left bottom}: intersection of $\mathcal{B}(0,\boldsymbol{\nu})$ and $\mathcal{B}(5,\boldsymbol{\nu})$ in green and $\mathcal{C}(0,\boldsymbol{\nu}) = \mathcal{B}(0,\boldsymbol{\nu})\setminus\mathcal{B}(5,\boldsymbol{\nu})$ in gray. Smoothing $M_0$ in the gray region is impossible without violating \eqref{eq: shift-cancel}. This leads to the distinction between regular (yellow) and singular (red) boundaries at $\xi_1 = \pm\pi/2$.\\\hspace{\textwidth}
 {\it Right}: boundary classification, singular (red) and regular (yellow) after similar arguments for all $\boldsymbol{\nu}$.}
\label{fig: boundary}
\end{figure}

However, because of the discontinuity of these shearlets across the boundary of their supports in the frequency domain, such shearlets have slow decay in the time domain. In order to improve their spatial localization, these boundaries need to be smoothed. In the remainder of this paper, we explore how this can be done for the first dyadic layer (without cutting further at higher frequencies); we call the corresponding basis {\it quasi-Shearlets}. As a result of smoothing the $M_j$'s, the shift cancellation no longer holds automatically as $supp(M_j)$ and $(supp(M_j)-\boldsymbol{\nu})$ may overlap in the neighborhood of the smoothed boundaries, see Fig. \ref{fig: boundary}, illustrating $M_j(\boldsymbol{\xi})\overline{M_j(\boldsymbol{\xi} + \boldsymbol{\nu})}\not\equiv 0.$

 Next, we show that there are regular boundaries which can be smoothed without violating shift cancellation, and singular boundaries which cannot be smoothed without violating constraints \eqref{eq: id-sum} and \eqref{eq: shift-cancel}.
 
\subsection{Boundary classification}
For simplicity, we introduce the following notations: let 
$\mathcal{B}(j,\boldsymbol{\nu}) =  supp(M_j)\cap\,(supp(M_j)-\boldsymbol{\nu})$ be the support of  $M_j(\boldsymbol{\xi})\overline{M_j(\boldsymbol{\xi} + \boldsymbol{\nu})}$ associated to $M_j$ and shift $\boldsymbol{\nu}$;
let $\mathcal{C}(j,\boldsymbol{\nu}) = \mathcal{B}(j,\boldsymbol{\nu}) \setminus \bigcup_{j'\neq j}\mathcal{B}(j',\boldsymbol{\nu})$. 
In addition, we may introduce the (slight abuse of) notation $\mathbf{C}_0 = \mathbf{S}_1$.

By definition, on $\mathcal{C}(j,\boldsymbol{\nu})$, $M_{j'}(\boldsymbol{\xi})\overline{M_{j'}(\boldsymbol{\xi}+\boldsymbol{\nu})} \equiv 0,\, \forall j' \neq j$. This implies the following lemma,

\begin{lem}
Shift cancellation can hold for $\boldsymbol{\nu}$, only if $\mathcal{C}(j,\boldsymbol{\nu})=\varnothing$. 
\end{lem}
Therefore, boundaries that after smoothing make $\mathcal{C}(j,\boldsymbol{\nu})$ non-empty are the {\it singular} boundaries, and the rest are {\it regular} boundaries. The next proposition provides an explicit method of classifying boundaries.
\begin{prop}\label{prop: class-bdy}
Let $supp(M_j) = \ov{$\mathbf{C}_j$}$,  then its boundary $\partial \mathbf{C}_j=\bigcup_{\boldsymbol{\lambda}\in \Lambda\setminus\{0\}}\mathcal{B}(j,\boldsymbol{\lambda})$. The set of singular boundaries of $\partial \mathbf{C}_j$ is $\bigcup_{\boldsymbol{\lambda}\in \Lambda\setminus\{0\}}\mathcal{C}(j,\boldsymbol{\lambda})$, whereas its compliment set is the regular boundary set.
 \end{prop}
 We use the notation $\mathcal{B}_s(j,\boldsymbol{\nu})$ for $\mathcal{B}(j,\boldsymbol{\nu})$ in the special case $supp(M_j) = \ov{$\mathbf{C}_j$}$ hereafter.
The boundary classification based on Proposition \ref{prop: class-bdy} is shown in the right of Fig. \ref{fig: boundary}, where the boundaries on the four corners of both $\mathbf{S}_0$ and $\mathbf{S}_1$ are singular: smoothing then is not allowed. On the other hand, the identity summation constraint \eqref{eq: id-sum} implies that $|M_j(\boldsymbol{\xi})| = 1$ on $\mathcal{C}(j,\boldsymbol{\nu}),\,\forall\boldsymbol{\nu}$, hence the discontinuity of $M$-functions across these boundaries is unavoidable.

\subsection{Pairwise smoothing of regular boundary}
Despite of the singular boundaries, better spatial localization can be achieved by carefully smoothing the regular boundaries. The regular boundaries of both $\mathbf{C}_{j_1}$ and $\mathbf{C}_{j_2}$ with adjacent supports consist of $\mathcal{B}_s(j_1,\boldsymbol{\nu})\cap\mathcal{B}_s(j_2,\boldsymbol{\nu})$, which can be denoted as a triple $(j_1,j_2,\boldsymbol{\nu})$. The following proposition shows that the regular boundaries $(j_1,j_2,\boldsymbol{\nu})$ can be paired according to shift pairs $(\boldsymbol{\nu},-\boldsymbol{\nu})$, and the boundaries must be smoothed pairwise within their $\epsilon-$neighborhood, $\mathcal{B}_{\epsilon}(j_1,j_2,\boldsymbol{\nu})$.

\begin{prop}\label{prop: pair-smooth}
Given $(j_1,j_2,\boldsymbol{\nu})\neq \emptyset$, $\exists \boldsymbol{\nu}',\,s.t. \boldsymbol{\nu}+\boldsymbol{\nu}'=0(mod(2\pi,2\pi))$ and $(j_1,j_2,\boldsymbol{\nu}')\neq\emptyset$. In addition, the identity summation and shift cancellation conditions hold if
\begin{itemize}
\item[{\it (1)}] $M_j = \mathbbm{1}_{\mathbf{C}_j},\,j\neq j_1,j_2$\\
\hspace*{-2em} and on $\mathcal{B}_{\epsilon}(j_1,j_2,\boldsymbol{\nu})\cup\mathcal{B}_{\epsilon}(j_1,j_2,\boldsymbol{\nu}')$\vspace*{.1em}
\item[{\it(2)}] $|M_{j_1}|^2 + |M_{j_2}|^2 = 1,$ 
\item[{\it (3)}] $\sum_{j_1,j_2} M_j(\cdot)\overline{M_j(\cdot+\boldsymbol{\nu})} = 0,$
\item[{\it (4)}] $\sum_{j_1,j_2} M_j(\cdot)\overline{M_j(\cdot+\boldsymbol{\nu}')} = 0.$ 
\end{itemize}
\end{prop}
Proposition \ref{prop: pair-smooth} enables us to smooth certain pairs of regular boundaries starting from the Shannon-type quasi-shearlets with the simplified conditions (2), (3) and (4); condition (1) can be removed as long as the initial $M_j$'s already satisfy identity summation and shift cancellation conditions and every $\boldsymbol{\xi}\in\mathbf{S}_0$ is not covered by more than two $M$ functions. Therefore, we can smooth regular boundaries pairwise, one after another.

The next proposition gives an explicit design of $(M_{j_1}, M_{j_2})$ satisfying the simplified conditions (2)-(4) in Proposition \ref{prop: pair-smooth} as well as a necessary condition of any valid design.
\begin{prop}\label{prop: M-design}
Given $(j_1,j_2,\boldsymbol{\nu})\neq\emptyset$, on $\Omega\cup(\Omega+\boldsymbol{\nu})$
\begin{itemize}
\item[(i)] $\sum_{j_1,j_2}M_{j}(\boldsymbol{\xi})\overline{M_{j}(\boldsymbol{\xi}+\boldsymbol{\nu})} = 0$
\item[(ii)] $\sum_{j_1,j_2}|M_{j}(\boldsymbol{\xi})|^2= 1 $ 
\end{itemize}
imply that\\
$|M_{j_1}(\boldsymbol{\xi})| = |M_{j_2}(\boldsymbol{\xi}+\boldsymbol{\nu})|,\quad |M_{j_2}(\boldsymbol{\xi})| = |M_{j_1}(\boldsymbol{\xi}+\boldsymbol{\nu})|$\\
Furthermore, if 
\[M_{j} = e^{i\boldsymbol{\xi}^{T}\boldsymbol{\eta}_{j}}\mathcal{M}_{j},\quad j=j_1,j_2, \text{ on }\Omega,\] where $\mathcal{M}_j$ is a real-valued function, 
$e^{i\boldsymbol{\nu}^T(\boldsymbol{\eta}_{j1}-\boldsymbol{\eta}_{j2})} = -1,$ and 
\[\mathcal{M}_{j_1}(\boldsymbol{\xi}) = \mathcal{M}_{j_2}(\boldsymbol{\xi}-\boldsymbol{\nu}),\;\mathcal{M}_{j_2}(\boldsymbol{\xi}) = \mathcal{M}_{j_1}(\boldsymbol{\xi}-\boldsymbol{\nu}),\text{ on }\Omega+\boldsymbol{\nu},\] 
then $(i)$ holds.
\end{prop}
Proposition \ref{prop: M-design} breaks down the design of $(M_{j_1},M_{j_2})$ into a pair of real function $(\mathcal{M}_{j_1}, \mathcal{M}_{j_2})$ on $\mathcal{B}_{\epsilon}(j_1,j_2,\boldsymbol{\nu})$ and a two phases $\boldsymbol{\eta}_1,\boldsymbol{\eta}_2$. 
The only constraint on $(\mathcal{M}_{j_1},\mathcal{M}_{j_2})$ for (ii) in Proposition \ref{prop: M-design} holds as well is that on $\mathcal{B}_{\epsilon}(j_1,j_2,\boldsymbol{\nu})$, $\sum_{j_1,j_2}|\mathcal{M}_{j}(\boldsymbol{\xi})|^2= 1$, which is easy to be satisfied.
We may construct all local pairs of $(\mathcal{M}_{j_1},\mathcal{M}_{j_2})$ separately, and put together afterwards different pieces of each $\mathcal{M}_j$ locate in different regular boundary neighborhoods $\mathcal{B}_{\epsilon}(j,j',\boldsymbol{\nu})$. 

\vspace*{.2em}
The phase term $e^{i\boldsymbol{\xi}^T\boldsymbol{\eta}_j}$ is preferrably defined on the full frequency domain, hence $\boldsymbol{\eta}_j$'s need to be solved globally. This global phase problem is stated precisely in the following proposition and one solution is provided.
\begin{prop}\label{prop: phase}
Applying Proposition \ref{prop: M-design} to all regular boundaries requires a set of phases $\{\boldsymbol{\eta}_j\}_{j = 0}^6,$ s.t.\[\textstyle e^{i\boldsymbol{\nu}^T(\boldsymbol{\eta}_{j1}-\boldsymbol{\eta}_{j2})} = -1, \quad\forall (j_1,j_2,\boldsymbol{\nu})\in\Delta,\]
{\small\begin{multline*}
\Delta = \Bigl\{\big(0,2,(0,\pi\big),\, \big(0,5,(\pi,0)\big),\,\big(1,3,(\pi,0)\big),\,\big(4,6,(0,\pi)\big),\\
\big(1,6,(\pi/2, 3\pi/2)\big),\,\big(2,3,(\pi/2, 3\pi/2)\big),\,\big(4,5,(\pi/2, 3\pi/2)\big), \\
\big(3,4,(\pi/2, \pi/2)\big),\big(1,2,(\pi/2, \pi/2)\big),\,\big(5,6,(\pi/2, \pi/2)\big)\Bigr\}
\end{multline*}}
There exsits a non-unique solution: 
{\small\begin{multline*}\boldsymbol{\eta}_0 = (0,0),\,\boldsymbol{\eta}_1 = (0,0),\,\boldsymbol{\eta}_2 = (1,1),\,\boldsymbol{\eta}_3 = (1,-1),\\
\boldsymbol{\eta}_4 = (0,2),\,\boldsymbol{\eta}_5=(1,1),\,\boldsymbol{\eta}_6 = (-1,1).
\end{multline*}}
\end{prop}

To summarize, Proposition \ref{prop: M-design} and \ref{prop: phase} introduce the following regular boundary smoothing scheme of $M$ function:
\begin{itemize}
\item[1.] Let $\mathcal{M}_j = \mathbbm{1}_{\mathbf{C}_j}$, and change its value to smooth a pair of regular boundaries $(j_1,j_2,\pm\boldsymbol{\nu})$ following step 2 and 3.
\item[2.]  On $\mathcal{B}_{\epsilon}(j_1,j_2,\boldsymbol{\nu})$,\\
\hspace*{2em} design $(\mathcal{M}_{j_1},\mathcal{M}_{j_2}),\quad$ s.t.
$\sum_{j_1,j_2}|\mathcal{M}_{j}(\boldsymbol{\xi})|^2= 1$.
\item[3.] On $\mathcal{B}_{\epsilon}(j_1,j_2,-\boldsymbol{\nu})$, \\
\hspace*{2em}let $\mathcal{M}_{j_1}(\boldsymbol{\xi}) = \mathcal{M}_{j_2}(\boldsymbol{\xi}-\boldsymbol{\nu})$, $\mathcal{M}_{j_2}(\boldsymbol{\xi}) = \mathcal{M}_{j_1}(\boldsymbol{\xi}-\boldsymbol{\nu})$\vspace*{.1em}
\item[4.] Repeat step 2 and 3 for all $(j_1,j_2,\boldsymbol{\nu})\in\Delta$. 
\item[5.]Let $M_j(\boldsymbol{\xi}) =e^{i\boldsymbol{\xi}^T\boldsymbol{\eta}_j} \mathcal{M}_j(\boldsymbol{\xi}),$ on $\mathbf{S}_0$, using the solution in Proposition \ref{prop: phase}.
\end{itemize}

\begin{figure}[!t]
\centering
\includegraphics[width=.15\textwidth]{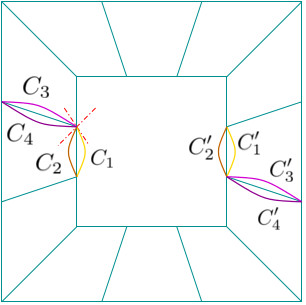}\hspace*{2mm}
\includegraphics[width=.25\textwidth]{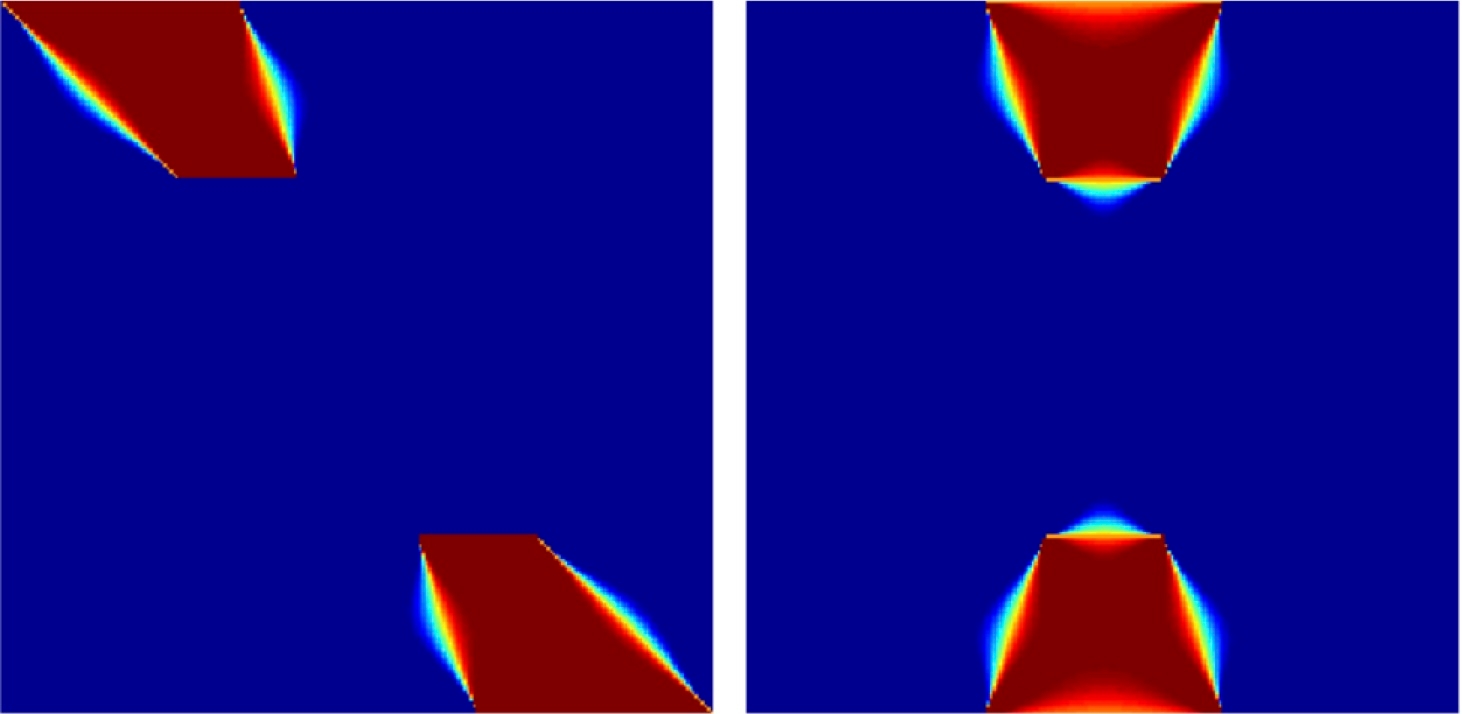}
\caption{ Left: contour design of $supp(M_5)$, Right: frequency support $|\hat{\psi}^j|$}
\label{fig: design}
\end{figure}

\section{Quasi-shearlet bases construction}
In this section, we present a family of quasi-shearlet bases constructed based on the $M$-function design discussed in section \ref{sec: design} using our proposed MRA framework.

We smooth all the regular boundaries except the ones on the boundary of $\mathbf{S}_0$. In the neighborhood of a regular boundary $\mathcal{B}_{\epsilon}(j,j',\boldsymbol{\nu})$, the change of $|M_j|$ from 0 to 1 depends on $\mathcal{M}_j$ and the contour of stop-band/pass-band is the boundary of level set $\{\mathcal{M}_j(\boldsymbol{\xi}) = 0\}\big/\{\mathcal{M}_j(\boldsymbol{\xi}) = 1\}$. Fig. \ref{fig: design} shows our design of the stop-band/pass-band contours of regular boundaries {\small $\big(5,6,(\pi/2,\pi/2)\big)$} and {\small $\big(0,5,(\pi,0)\big)$}. It is necessary that the contours intersect only at the vertices of $\mathbf{C}_5$, so that $supp(M_5)\cap supp(M_6)\cap supp(M_0) = \emptyset$ as the relaxed condition (1) in Proposition \ref{prop: pair-smooth}. Moreover, $\mathcal{M}_5$ is designed to be symmetric with respect to the origin in both regular boundary neighborhoods. 

\begin{figure}[!t]
\centering
\includegraphics[width=.3\textwidth]{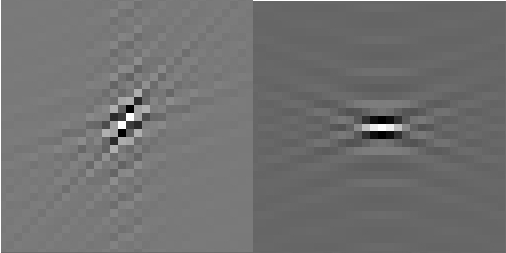}\hspace*{2mm}
\includegraphics[width=.1\textwidth]{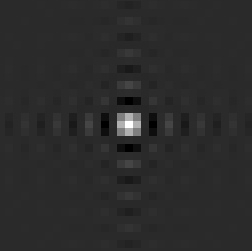}
\caption{ Left: quasi-shearlet $\psi^j$, Right: scaling function $\phi$}
\label{fig: quasi-shear}
\end{figure}

The contours related to other regular boundaries are designed likewise to achieve the best symmetry, and in particular, the quasi-shearlets are real. The right of Fig.\ref{fig: design} shows the frequecy support of quasi-shearlets generated by such design; Fig.\ref{fig: quasi-shear} shows the quasi-shearlets and scaling function in space domain. This Parseval frame is indeed an orthonormal basis by checking Theorem \ref{thm: basis cond}.

As shown in Fig. \ref{fig: quasi-shear}, the quasi-shearlets orient in six directions, yet they are not spatially localized due to two limitations of our framework of design. The major limitation is the singular boundaries on the corners of the low-frequency square $\mathbf{S}_1$, where the discontinuity in frequency domain is inevitable. The minor limitation is the lack of smoothness at the vertices of $M_2$ and $M_5$ because of the assumption of no triple overlapping of $M$ functions. This could possibly be avoided by using a more delicate (but more complicated) design around the vertices $(\pm\frac{\pi}{2},\pm\frac{\pi}{6})$.

\begin{figure}[!t]
\centering
\includegraphics[width=.3\textwidth]{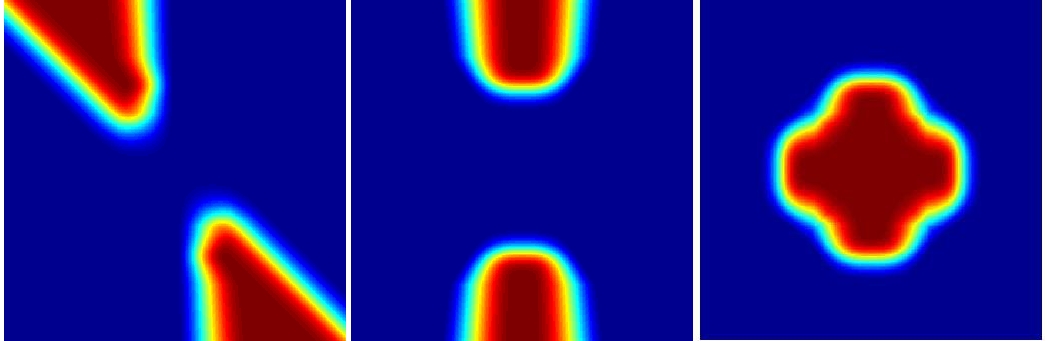}\hspace*{2mm}
\includegraphics[width=.1\textwidth]{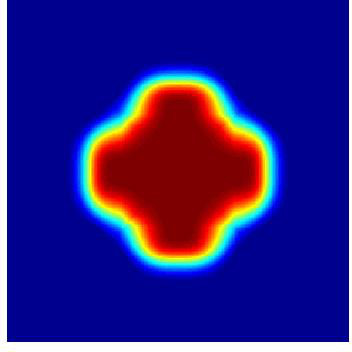}
\caption{ $M$ function of a quasi-shearlet tight frame with dyadic downsampling scheme. Left: $M_j,\,j=1,\cdots,6$, Right: $M_0$}
\label{fig: tightframe-M}
\end{figure}

\begin{figure}[!t]
\centering
\includegraphics[width=.3\textwidth]{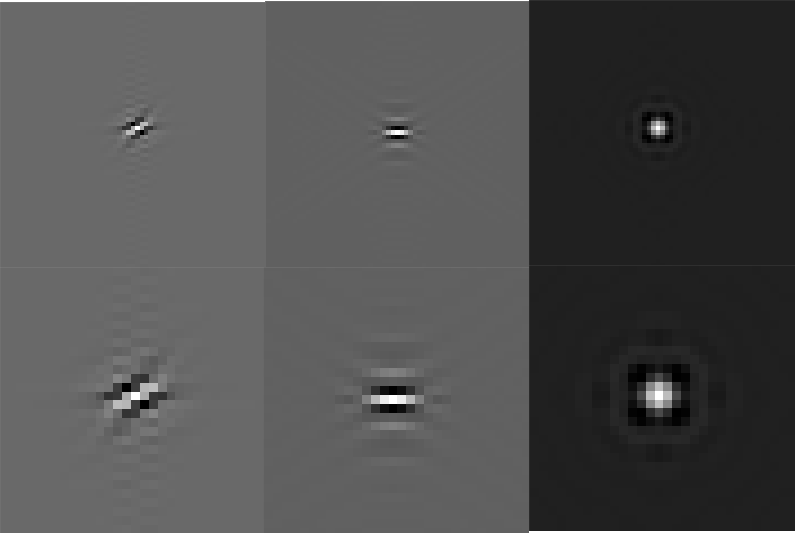}
\caption{ $\psi^1,\,\psi^2$ and $\phi$ in a quasi-shearlet tight frame with dyadic downsampling scheme. Top: same resolution as in Fig. \ref{fig: quasi-shear}. Bottom: higher resolution of the top row. }
\label{fig: tightframe-sp}
\end{figure}

On the other hand, if we allow a bit redundancy and abandon the critical downsampling scheme, then the corresponding set of shifts $\Lambda$ will be smaller, thus fewer shift cancellation constraints \eqref{eq: shift-cancel} need to be satisfied. Therefore, the discontinuities of a quasi-shearlets basis in the frequency domain around the singular boundaries can be removed in a low redundant quasi-shearlet tight frame. In particular, we construct a quasi-shearlet tight frame with redundancy of 2 by using the classical dyadic downsampling $D_2$, with shift cancellaiton constraints \eqref{eq: shift-cancel} only on set $\Gamma\setminus\{\boldsymbol{0}\}$. Even though $supp(M_0)$ still cannot be extended outside of the four corners of $\mathbf{C}_0$ due to \eqref{eq: shift-cancel}, $M_1,M_3,M_4$ and $M_6$ can cover the inside of these corners now. This makes smoothing the boundaries of $M_0$ inwards possible without violating \eqref{eq: id-sum}, see Fig. \ref{fig: tightframe-M}. With a price of double redundancy, we obtain shearlets with much better spatial localization as shown in Fig. \ref{fig: tightframe-sp}.

\section{Conclusion and future work}
We observed that it is possible to construct a Shannon-type orthonormal shearlet basis and characterized the corresponding subsampling scheme,  
which involves quincunx downsampling following the dyadic downsampling.
We then investigated whether these can be smoothed in the frequency domain, in order to obtain better spatial localization. Restricting our discussion to the first dyadic frequency ring (and thus to what we have called quasi-shearlets), we showed that some smoothing is possible, but some discontinuities in the frequency domain must persist.

This scheme can be relaxed to provide shearlets with much better spatial localization and low redundancy.


\bibliographystyle{IEEEtran}
\bibliography{shearlet}

\end{document}